\documentclass[12pt]{article}

\usepackage{amsfonts}

\usepackage{amsmath,amssymb}

 \numberwithin{equation}{subsection}

\begin{document}

\title{Some multiplication formulas in an affine Hecke algebra   }
\author{Liping Wang\\
\small China Economics and Management Academy,\\
\small Central University of Finance and Economics, Beijing, China\\
\small wanglp@amss.ac.cn}

\date{}

\maketitle

\begin{abstract}

In this paper we consider the Hecke algebra $\mathcal {H}$
associated to an extended affine Weyl group of type
$\widetilde{B_2}$. We give some interesting formulas on
$C_{rt}S_{\lambda}$, which imply some relations between the
Kazhdan-Lusztig coefficients $\mu (y,w)$ and representations of some
algebraic groups. Here $C_{rt}$ is one element in the
Kazhdan-Lusztig basis of $\mathcal {H}$ and $S_{\lambda}$ is an
element in the center of $\mathcal {H}$.

\end{abstract}

\maketitle
\section*{Introduction} In [W], the author studies the leading
coefficients of the Kazhdan-Lusztig polynomials for
  an Affine Weyl group of type $\widetilde{B_2}$. But in that paper
  not all of the leading coefficients are given. In order to
  complete the problem, we study some multiplication formulas in the
  associated Hecke algebra in this paper. By these formulas, we will
  see some close relations between the
Kazhdan-Lusztig coefficients $\mu (y,w)$ and representations of some
algebraic groups.
\section{Preliminaries} In this section we recall some basic facts about Hecke algebras
which will be needed later.
\subsection{Hecke algebras} Let $G$ be a connected, simply connected reductive
algebraic group over the field $\mathbb{C}$ of complex numbers and
$T$ a maximal torus of $G.$ Let $N_G (T)$ be the normalizer of $T$
in $G.$ Then $W_0 =N_G (T)/T$ is a Weyl group, which acts on the
character group $\Lambda=\textrm{Hom}(T, \mathbb{C} ^*)$ of $T$. Let
$\Lambda_r$ be the root lattic of $G,$ then the semi-direct product
$W'=W_0 \ltimes \Lambda_r$ is an affine Wely group and $W=W_0
\ltimes \Lambda$ is called an extended affine Weyl group associated
with $W'$. $W'$ is a Coxeter group, while $W$ is not in general.

We shall denote by $S$ the set of simple reflections of $W'.$ We can
find an abelian subgroup $\Omega$ of $W$ such that $\omega
S=S\omega$ for any $\omega \in\Omega$ and $W=\Omega \ltimes W'.$ We
shall denote the length function of $W'$ by $l$ and use $\leq$ for
the Bruhat order on $W'.$ The length function $l$ and the partial
order $\leq$ on $W'$ can be extended to $W$ as usual, that is,
$l(\omega w)=l(w),$ and $\omega w\leq \omega' u$ if and only if
$\omega=\omega'$ and $w\leq u,$ where $\omega, \omega'$ are in
$\Omega$ and $w,u$ are in $W'.$

Let $\cal H'$ be the Hecke algebra of $(W',S)$ over $\mathcal
{A}=\mathbb{Z}$ $[q^{\frac{1}{2}}, q^{-\frac{1}{2}}]$ ($q$ an
indeterminate) with parameter $q.$ Let $\{T_w\}_{w\in W'}$ be its
standard basis and $C_w =q^{-\frac{l(w)}{2}}\sum_{y\leq w}
P_{y,w}T_y, w\in W$ be its Kazhdan-Lusztig basis, where $P_{y,w}\in
\mathbb{Z} [q]$ are the Kazhdan-Lusztig polynomials. The degree of
$P_{y,w}$ is less than or equal to $\frac{1}{2} (l(w)-l(y)-1)$ if
$y<w$ and $P_{w,w}=1.$

For an affine Weyl group $W'$ , we know that the coefficients of
these polynomials are all non-negative (see [KL2]).

We write $P_{y,w}=\mu (y,w)q^{\frac{1}{2} (l(w)-l(y)-1)}+$ lower
degree terms. The coefficient $\mu (y,w)$ is very interesting, this
can be seen even from the recursive formula for Kazhdan-Lusztig
polynomials ( see [KL1]). We call $\mu (y,w)$ the Kazhdan-Lusztig
 coefficient of $P_{y,w}.$ We denote by $y\prec w$ if
$y\leq w$ and $\mu (y,w)\neq0.$ Define $\widetilde{\mu }(y,w)=\mu
(y,w)$ if $y\leq w$ or $\widetilde{\mu }(y,w)=\mu (w,y)$ if $w\leq
y.$

Let $\cal H$ be the generic Hecke algebra of $W.$ Then the algebra
$\cal H$ is isomorphic to the ``twisted" tensor product $\mathbb{Z}
[\Omega]\otimes_{\mathbb{Z}}\cal H'$.

We have the following multiplication formula (see [KL1]):\\
Given any element $w\in W$, then\\
(a) For $s\in S$ we have
$${C_sC_w}=\left\{\begin{array}{ll}
{(q^{\frac{1}{2}}+q^{\frac{1}{2}})C_w},& \ \ \textrm{if} \;  sw<w ,\\
{C_{sw}+\sum_{sy<y\prec w }\mu (y,w)C_y},& \ \ \textrm{if}\ sw>w
\end{array} \right.$$
and$${C_wC_s}=\left\{\begin{array}{ll}
{(q^{\frac{1}{2}}+q^{\frac{1}{2}})C_w},& \ \ \textrm{if} \;  ws<w ,\\
{C_{ws}+\sum_{ys<y\prec w }\mu (y,w)C_y},& \ \ \textrm{if}\ ws>w. \end{array} \right.$$\\

We refer to [KL1] for the definition of the preorders $\leq_L,
\leq_R, \leq_{LR}$ and of the equivalence relations $\sim_L, \sim_R,
\sim_{LR}$ on $W'.$ The corresponding equivalence classes are called
left cells, right cells, two-sided cells of $W',$ respectively.

For any $u=\omega_1u_1,\ w=\omega_2w_1,\ \omega_1,\omega_2\in
\Omega,\ u_1,w_1\in {W'},$ we define $P_{u,w}=P_{u_1,w_1}$ if
$\omega_1=\omega_2$ and define $P_{u,w}=0$ if
$\omega_1\neq\omega_2.$ We say that $u\leq_Lw$ or $u\leq_Rw$ or
$u\leq_{LR}w$ if $u_1\leq_Lw_1$ or
$\omega_1u_1\omega_1^{-1}\leq_R\omega_2w_1\omega_2^{-1}$ or
$u_1\leq_{LR}w_1.$ The left (resp. right or two-sided) cells of $W$
are defined as those of $W'$. We also define $a(\omega w)=a(w)$ for
$\omega\in\Omega,
 w\in W'$.

For $w\in W$, set $L(w)=\{s\in S\mid sw\leq w\},\;R(w)=\{s\in S\mid
ws\leq w\}.$ Then we have (see [KL1])\\
(b) $R(w)\subseteq R(y), \;\textrm{if}\;y\leq_L w.$ In particular,
$R(w)=R(y),\;\textrm{if}\;y\sim_L w;$\\
(c) $L(w)\subseteq L(y), \;\textrm{if}\;y\leq_R w.$ In particular,
$L(w)=L(y),\;\textrm{if}\;y\sim_R w$.

\noindent {\bf Convention:}   For any element in $\Lambda$, we will
use the same notation when it is regarded as an element in $W$. And
we will use addition and multiplication in $\Lambda$ and $W$
respectively.
\subsection{The center of a Hecke algebra} In [L3], Bernstein described the center of Hecke algebras. We
recall it in the following.

To each $x\in \Lambda$, Bernstein associates an element
$\theta_{x}\in \mathcal {H}$ defined by
$$\theta_{x}=(q^{-l(x')/2}T_{x'})(q^{-l(x'')/2}T_{x''})^{-1}, $$where
$x',x''$ are elements in $\Lambda^+$ such that $x=x'-x''$ ( or we
write $x=x'x''^{-1}$ in $\widetilde{W}$), where $\Lambda^+$ is the
dominant weights set which is the set $\{x\in\Lambda\mid
l(xw_0)=l(x)+l(w_0)\}.$ $\theta_{x}$ is well defined and  is
independent of the choice of $x',x''$ (for details see [L3]).

Let $x\in\Lambda^+$ and $W_0(x)$ be its $W_0-$orbit in $\Lambda$. We
denote $z_{x}=\sum_{x'\in W_0(x)}\theta_{x'}$.

For $x\in\Lambda^+$, we define $$S_x=\sum_{{x'\in \Lambda^+}\atop
x'\leq x}d_{x'}(x)z_{x'},$$ where $d_{x'}(x)$ means the multiplicity
of weight $x'$ in $V(x)$, and $V(x)$ is a rational irreducible
$G-$module of highest weight $x$.

Bernstein showed that $S_x,x\in\Lambda^+$ form an $\mathcal
{A}$-basis of the center of $\mathcal {H}$.

Then we have (see [L3])\\
(a) $S_xS_{x'}=\sum_{z\in\Lambda^+}m_{x,x',z}S_z$ for any
$x,x'\in\Lambda^+$. Here $m_{x,x',z}$ is defined to be  the
multiplicity of $V(z)$ in the tensor product $V(x)\otimes V(x').$
\section{Cells in an affine Weyl group of type $\widetilde{B_2}$}

 In the rest of this
paper,  $G=\textrm{Sp}_4$$(\mathbb{C})$. Then $(W',S)$ is an affine
Weyl group of type $\widetilde{B_2}$ and $W$ is the extended affine
Weyl group associated with $\textrm{Sp}_4$$(\mathbb{C})$. Let
$S=\{r, s, t\}$ be the set of all simple reflections of $W'$ with
$rt=tr,\ (rs)^4 =(st)^4 =e.$ The Weyl group $W_0$ is generated by
$s$ and $t$.

We assume that $s$ is the simple reflection of $W_0$ corresponding
to the long simple root $\alpha_1$, and $t$ is the simple reflection
corresponding to the short simple root $\alpha_2$.

We have that $x_1=\alpha_1+\alpha_2$ and
$x_2=\frac{1}{2}\alpha_1+\alpha_2$ are the fundamental dominant
weights, the corresponding elements in $W$ are $x_1=stsr$ and
$x_2=\omega rsr$ respectively, where $\omega\in \Omega$ such that
$W=\Omega\ltimes W',\ \Omega=\{e,\omega\},$ $\omega
r=t\omega,\;\omega s=s\omega\ \textrm{and}\ \omega t=r\omega$.

In [L1], Lusztig described the left cells and two-sided cells of
$(W',S).$ For any subset $J$ of $S=\{r, s, t\}$, we denote by $W^J$
the set of all $w\in W'$ such that $R(w)=J.$ Then $(W',S)$ has $16$
left cells:
$$A_{rs}=W^{rs},\;A_{rt}=A_{rs}t,\;A_s=A_{rt}s,\;A_r=A_sr,$$
$$A_{st}=W^{st},\;A_{rt}'=A_{st}r,\;A_s'=A_{rt}'s,\;A_t=A_s't,$$
$$B_{rt}=W^{rt}-(A_{rt}\cup A_{rt}'),\;B_s=B_{rt}s,\;B_r=B_sr,$$
$$\;\;\;\;\;\;\;\;\;\;\;B_t=B_st,\;C_r=W^r-(A_r\cup B_r),\;C_t=W^t-(A_t\cup B_t),$$
$$C_s=W^s-(A_s\cup A_s'\cup B_s),\;D_\emptyset=W^\emptyset=\{e\}.\;\;\;\;\;\;\;\;$$

Set $$c_e=D_\emptyset,\ c_1=C_r\cup C_s\cup C_t,\ c_2=B_r\cup
B_s\cup B_t\cup B_{rt},$$
$$c_0=A_r\cup
A_s\cup A_s'\cup A_t\cup A_{rs}\cup A_{st}\cup A_{rt}\cup A_{rt}'.$$

From [L1], we know that $c_e,\ c_1,\ c_2,\ c_0$ exhaust two-sided
cells of $W'.$ We have that
\begin{eqnarray*}c_e
&=&\{x\in W'\mid a(x)=0\}=\{e\}
\\
c_1&=&\{x\in W'\mid a(x)=1\}\\
c_2&=&\{x\in W'\mid a(x)=2\}\\
c_0&=&\{x\in W'\mid a(x)=4\}.
\end{eqnarray*}

So the two-sided cells of $W$ exactly are
\begin{eqnarray*}c_e
&=&\{x\in W\mid a(x)=0\}=\{e\}
\\
c_1&=&\{x\in W\mid a(x)=1\}\\
c_2&=&\{x\in W\mid a(x)=2\}\\
c_0&=&\{x\in W\mid a(x)=4\}.
\end{eqnarray*}


\section {Main results}  \ \ \ \ \ In this section, we
consider the coefficients of the product $C_{rt}S_{\lambda}$ under
the basis $\{C_w|\ w\in W\}$ for any dominant weight $\lambda$,
where $S_{\lambda}$ is an element in the center of $\mathcal {H}$
(see the definition of $S_{\lambda}$ in Section 1.2). By the formula
for $C_{rt}S_{\lambda}$ in Theorem 3.6, we can see that the leading
coefficients of the Kazhdan-Lusztig polynomials have close relations
with representations of algebraic groups.

First we give some observations in the Hecke algebra $\mathcal {H}$
associated with the extended affine Weyl group $W$ of type
$\widetilde{B_2}$.

Let $\mathcal {H}_{c_0}=\sum_{w\in c_0}\mathcal {A}C_w$, where $c_0$
is the lowest two-sided cell of $W$. By [L1], we know that $\mathcal
{H}_{c_0}$ is a two-sided ideal of $\mathcal {H}$. We set
$[2]=q^{\frac{1}{2}}+q^{-\frac{1}{2}}$.

Then we have the following result.

\noindent {\bf Lemma 3.1} \\
(a) For any $1\leq m\in \mathbb{N}$, we
have\\
$$C_{rtsrt}C_{rt(srt)^m}=[2]^2(C_{rt(srt)^{m+1}}+C_{rt(srt)^{m-1}})\ \textrm{mod}\ \mathcal {H}_{c_0}.$$
(b) For any $1\leq m, n\in \mathbb{N}$, we
have\\
$$C_{rt(srt)^m}C_{rt(srt)^n}=[2]^2\sum_{i=0}^{\textrm{min}\{m,n\}}C_{rt(srt)^{m+n-2i}}\
\textrm{mod}\ \mathcal {H}_{c_0}.$$

\noindent {\bf Proof.} We prove (a) first. For any $1\leq m\in
\mathbb{N}$,  by the results in [W, Section 4.2, 4.3] and the
formula in 1.1(a), we have
\begin{eqnarray*}
& &C_{rtsrt}C_{rt(srt)^m}\\
&=&(C_{rts}-[2])C_{rt}C_{rt(srt)^m}\\
&=&[2]^2C_{rt}C_{(srt)^{m+1}}-[2]^3C_{rt(srt)^m}\\
&=&[2]^2C_{r}C_{t}C_{(srt)^{m+1}}-[2]^3C_{rt(srt)^m}\\
&=&[2]^2C_{r}(C_{t(srt)^{m+1}}+C_{rt(srt)^m})-[2]^3C_{rt(srt)^m}\ \textrm{mod}\ \mathcal {H}_{c_0}\\
&=&[2]^2C_{r}C_{t(srt)^{m+1}}\ \textrm{mod}\ \mathcal {H}_{c_0} \\
&=&[2]^2(C_{rt(srt)^{m+1}}+C_{rt(srt)^{m-1}})\ \textrm{mod}\
\mathcal {H}_{c_0}
\end{eqnarray*}

We can prove (b) by induction on $m$ for any $n\geq1$, with
assumption that $m\leq n$. The computation is very similar to the
above. \hfill$\Box$\

By the definition of $S_{\lambda}$ in Subsection 1.2 and some
computations, we have the following result.

\noindent {\bf Lemma 3.2} For the fundamental dominant weights $x_1$
and $x_2$, we have that
$$C_{rt}S_{x_1}=C_{rststr}+C_{tsrsrt}-[2]C_{rtsrt}+C_{rt}$$
and
$$C_{rt}S_{x_2}=(C_{rtsrt}-[2]C_{rt})C_{\omega}.$$

Let $\mathcal {H}_{\geq2}=\sum_{{w\in W}\atop a(w)\geq2}\mathcal
{A}C_w$. Obviously, $\mathcal {H}_{c_0}\subseteq\mathcal
{H}_{\geq2}$. By [L1], we know that $\mathcal {H}_{\geq2}$ is a
two-sided ideal of $\mathcal {H}$. Then $\mathcal
{H}_{\geq2}/\mathcal {H}_{c_0}$ has a natural $\mathcal {H}$-mod
structure with $\mathcal {A}$-basis $\{\widehat{C}_w|\ w\in c_2\}$,
where $\widehat{C}_w$ is the image of $C_w$ in $\mathcal
{H}_{\geq2}/\mathcal {H}_{c_0}$.

Let $\mathcal {H}_1$ be the $\mathcal {A}$-submodule of $\mathcal
{H}_{\geq2}/\mathcal {H}_{c_0}$ spanned by the elements\\
$\{\frac{1}{[2]^2}\widehat{C}_{rt(srt)^m\omega^p}|\ p=0,1;\ m\in
\mathbb{N}\}$. By Lemma 3.1(b), we know that $\mathcal {H}_1$ has an
algebra structure.

Let $F_{c_2}=SL_2(\mathbb{C})\times \mathbb{Z}/(2)$. We know that
the sets of irreducible representations up to isomorphism of
$SL_2(\mathbb{C})$ and $\mathbb{Z}/(2)$ are
$\textrm{Irr}_{SL_2(\mathbb{C})}=\{V(m)|\ m\in \mathbb{N}\}$ and
$\textrm{Irr}_{\mathbb{Z}/(2)}=\{1,\epsilon|\ \epsilon^2=1\}$,
respectively. Here $V(m)$ is an irreducible representation of
$SL_2(\mathbb{C})$ with highest weight $m$ and $\epsilon$ is the
 sign representation of $\mathbb{Z}/(2)$.
  Then the set of irreducible representations up to isomorphism of $F_{c_2}$ is
$\textrm{Irr}_{F_{c_2}}=\textrm{Irr}_{SL_2(\mathbb{C})}\times
\textrm{Irr}_{\mathbb{Z}/(2)}$.

Let $\textrm{R}_{F_{c_2}},\textrm{R}_{SL_2(\mathbb{C})}
\textrm{and}\ \textrm{R}_{\mathbb{Z}/(2)}$ be the rational
representation rings of $F_{c_2}$, $SL_2(\mathbb{C})$ and
$\mathbb{Z}/(2)$ respectively. We know that they are generated by
$\textrm{Irr}_{F_{c_2}}$, $\textrm{Irr}_{SL_2(\mathbb{C})}$ and
$\textrm{Irr}_{\mathbb{Z}/(2)}$, respectively. Then we have that
$\textrm{R}_{F_{c_2}}=\textrm{R}_{SL_2(\mathbb{C})}\times
\textrm{R}_{\mathbb{Z}/(2)}$

\noindent {\bf Proposition 3.3} The map $\varphi_1:\ \mathcal
{H}_1\longrightarrow {\mathcal A}\otimes R_{F_{c_2}}$ defined by\\
$\varphi_1(\frac{1}{[2]^2}\widehat{C}_{rt(srt)^m\omega^p})=V(m)\epsilon^p$
 is an isomorphism of $\mathcal {A}$-algebras.

\noindent {\bf Proof.} Obviously, the map $\varphi_1$ is a
bijection. It's enough to prove that $\varphi_1$ is a homomorphism
of $\mathcal {A}$-algebras. By Lemma 3.1(b), we have that
\begin{eqnarray*}
&&(\frac{1}{[2]^2}\widehat{C}_{rt(srt)^m\omega^p})(\frac{1}{[2]^2}\widehat{C}_{rt(srt)^n\omega^{p'}})\\
&=&\frac{1}{[2]^4}\widehat{C_{rt(srt)^m\omega^p}C_{rt(srt)^n\omega^{p'}}}\\
&=&\frac{1}{[2]^2}\widehat{C}_{\omega^p\omega^{p'}}\sum_{i=0}^{\textrm{min}\{n,m\}}
\widehat{C}_{rt(srt)^{m+n-2i}}.
\end{eqnarray*}
Thus
$$\varphi_1((\frac{1}{[2]^2}\widehat{C}_{rt(srt)^m\omega^p})(\frac{1}{[2]^2}\widehat{C}_{rt(srt)^n\omega^{p'}}))=
\epsilon^p\epsilon^{p'}\sum_{i=0}^{\textrm{min}\{n,m\}}V(m+n-2i).$$
On the other hand,
\begin{eqnarray*}
&&\varphi_1(\frac{1}{[2]^2}\widehat{C}_{rt(srt)^m\omega^p})\varphi_1(\frac{1}{[2]^2}\widehat{C}_{rt(srt)^n\omega^{p'}})\\
&=&\epsilon^p\epsilon^{p'}V(m)\otimes V(n)\\
&=&\epsilon^p\epsilon^{p'}\sum_{i=0}^{\textrm{min}\{n,m\}}V(m+n-2i).
\end{eqnarray*}
We get that
$$\varphi_1(\frac{1}{[2]^2}\widehat{C}_{rt(srt)^m\omega^p})\varphi_1(\frac{1}{[2]^2}\widehat{C}_{rt(srt)^n\omega^{p'}})
=\varphi_1((\frac{1}{[2]^2}\widehat{C}_{rt(srt)^m\omega^p})(\frac{1}{[2]^2}\widehat{C}_{rt(srt)^n\omega^{p'}}))$$
Thus $\varphi_1$ is an isomorphism of $\mathcal {A}$-algebras.
\hfill$\Box$\

Let $\mathcal {H}_2$ be the $\mathcal {A}$-submodule of $\mathcal
{H}_{\geq2}/\mathcal {H}_{c_0}$ spanned by the elements\\
$\{\frac{1}{[2]^2}\widehat{C_{rt}S_{\lambda}}|\ \lambda\in
\Lambda^+\}$, where $\widehat{C_{rt}S_{\lambda}}$ is the image of
$C_{rt}S_{\lambda}$ in $\mathcal {H}_{\geq2}/\mathcal {H}_{c_0}$.
 And by the fact that (see Section 1.2(a))
 $$(\frac{1}{[2]^2}C_{rt}S_{\lambda})(\frac{1}{[2]^2}C_{rt}S_{\lambda'})
 =\frac{1}{[2]^2}C_{rt}\sum_{z\in\Lambda^+}m_{\lambda,\lambda',z}S_z$$
we know that $\mathcal {H}_2$ is also an $\mathcal {A}$-algebra.

Obviously, we get the following result.

\noindent {\bf Proposition 3.4} The map $\varphi_2:\ \mathcal
{H}_2\longrightarrow \mathcal {A}\otimes R_{G}$ defined by
$\varphi_2(\frac{1}{[2]^2}\widehat{C_{rt}S_{\lambda}})=S_{\lambda}$
 is an isomorphism of $\mathcal {A}$-algebras, where
 $G$=$\textrm{Sp}_4$$(\mathbb{C})$.

\noindent {\bf Proof.} The only thing to note is the fact that we
can regard $S_{\lambda}$ as the character of the irreducible
$G-$module with highest weight $\lambda\in\Lambda^+$.\hfill$\Box$\

We also have  the following fact.

\noindent {\bf Lemma 3.5} We have that $\mathcal {H}_2\subseteq
\mathcal {H}_1$.

\noindent {\bf Proof.} For any $\lambda\in \Lambda^+$, we know that
$S_{\lambda}$ belongs to the center of $\mathcal {H}$. We can assume
that $C_{rt}S_{\lambda}=S_{\lambda}C_{rt}=\sum_{w\in
\widetilde{W}}a_wC_w$, where $a_w\in \mathcal {A}$. By [L1], we know
that $a_w\neq0$ implies that $w\leq_Rrt$ and $w\leq_Lrt$. Thus
$C_{rt}S_{\lambda}\in \mathcal {H}_{\geq2}$. And by 1.1 (b) and (c),
we get that $\{r,t\}\subseteq L(w)\cap R(w)$. Following the
decomposition of left cells of $W$ in Section 2, we get that
$a_w\neq0$ implies that $L(w)= R(w)=\{r,t\}$. Thus for any
$\lambda\in \Lambda^+$,
$\frac{1}{[2]^2}\widehat{C_{rt}S_{\lambda}}\in \mathcal {H}_1$.
\hfill$\Box$\

Define $\varphi=\varphi_1\varphi_2^{-1}:\ \mathcal {A}\otimes
R_{G}\longrightarrow \mathcal {A}\otimes R_{F_{c_2}}$. It's an
injective homomorphism of $\mathcal {A}$-algebras.

For any $\lambda\in \Lambda^+$, since $\varphi (S_{\lambda})\in
\mathcal {A}\otimes R_{F_{c_2}}$, thus we can assume that $\varphi
(S_{\lambda})=\sum_{{p=0,1}\atop k\in
\mathbb{N}}a_{k,p}V(k)\epsilon^p$, for some $a_{k,p}\in\mathcal
{A}$. In fact, for any given $\lambda\in \Lambda^+$, we have two
different ways to solve these $a_{k,p}$.

On the one side, by Lemma 3.2 we get that
\begin{eqnarray*}
&&\varphi(S_{x_1})=\varphi_1\varphi_2^{-1}(S_{x_1})=\varphi_1(\frac{1}{[2]^2}\widehat{C_{rt}S}_{x_1})\\
&=&
\varphi_1(\frac{1}{[2]^2}(\widehat{C}_{rststr}+\widehat{C}_{tsrsrt}-[2]\widehat{C}_{rtsrt}+\widehat{C}_{rt}))\\
&=& -[2]V(1)+1
\end{eqnarray*}
 and
\begin{eqnarray*}
&&\varphi(S_{x_2})=\varphi_1\varphi_2^{-1}(S_{x_2})=\varphi_1(\frac{1}{[2]^2}\widehat{C_{rt}S}_{x_2})\\
&=&
\varphi_1(\frac{1}{[2]^2}(C_{rtsrt}-[2]C_{rt})C_{\omega})\\
&=& (V(1)-[2])\epsilon.
\end{eqnarray*}
For any $\lambda\in \Lambda^+$, we know that
$S_{\lambda}=f(S_{x_1},S_{x_2})\in\mathbb{Z}[S_{x_1},S_{x_2}]$, the
polynomial ring of $S_{x_1}$ and $S_{x_2}$. Thus we have that

$$\varphi(S_{\lambda})=f(\varphi(S_{x_1}),\varphi(S_{x_2}))=f(-[2]V(1)+1,(V(1)-[2])\epsilon).$$
In $R_{SL_2(\mathbb{C})}$, we have that
$$V(1)^k\in\mathbb{Z}[V(1),\ldots,V(k)],\ \ \textrm{for\ any}\ k\in\mathbb{N}.$$
Thus $\varphi (S_{\lambda})=\sum_{{p=0,1}\atop k\in
\mathbb{N}}a_{k,p}V(k)\epsilon^p$, for some $a_{k,p}\in\mathcal
{A}$. Then we can get those $a_{k,p}$. This method may be not very
easy, because the polynomial $f(S_{x_1},S_{x_2})$ is not easy to
compute for general $S_{\lambda}$.

On the other side, we consider the question in $\mathcal {A}\otimes
R_{F_{c_2}}$. In fact we know that $F_{c_2}$ is a maximal reductive
subgroup of $C_{G}(u)$ for some unipotent element $u$ in $G$ which
corresponds to $c_2$. We denote by $T$ and $T'$ the maximal torus of
$G$ and $SL_2(\mathbb{C})$, respectively. Let $\Lambda=X(T)$ and
$\Lambda'=X(T')$ be the character groups. We have that
$\Lambda=\mathbb{Z}x_1+\mathbb{Z}x_2$ and $\Lambda'=\mathbb{Z}\xi$,
where $x_1$ and $x_2$ map diag$(\ a,b,a^{-1},b^{-1})$ to $ab$ and
$b$, respectively. We also have that $\xi$ is the character maps
diag$(\ a,a^{-1})$ to $a$. We have two group algebras
$\mathbb{Z}[\Lambda]=\textrm{Span}_{\mathbb{Z}}\{\theta_{x}|x\in\Lambda\}$
 and
 $\mathbb{Z}[\Lambda']=\textrm{Span}_{\mathbb{Z}}\{\theta_{x'}'|x'\in\Lambda'\}$.
We know that $\mathbb{Z}[\Lambda]$ is generated by $\theta_{x_1}$
and $\theta_{x_2}$, while $\mathbb{Z}[\Lambda']$ is generated by
$\theta_{\xi}'$. We can identify $R_G$ with
$\mathbb{Z}[\Lambda]^{W_0}$ and identify $R_{F_{c_2}}$ with
$\mathbb{Z}[\Lambda']^{S_2}\times\{e,\epsilon\}$.

We have a homomorphism of $\mathcal {A}-$algebras
$\widetilde{\varphi}:\ \mathcal {A}\otimes
\mathbb{Z}[\Lambda]\rightarrow \mathcal {A}\otimes
\mathbb{Z}[\Lambda']\times\{e,\epsilon\}$, which restricts to
$\mathcal {A}\otimes R_G$ is just $\varphi$. We can get that
$\widetilde{\varphi}(\theta_{x_1})=-q^{\frac{1}{2}}\theta_{\xi}'$
and $\widetilde{\varphi}(\theta_{x_2})=\theta_{\xi}'\times\epsilon$,
since we have gotten that
$\widetilde{\varphi}(S_{x_1})=-[2]V(\xi)+1$ and
$\widetilde{\varphi}(S_{x_2})=(V(\xi)-[2])\epsilon$.

For any $\lambda\in \Lambda^+$, we know that
$S_{\lambda}=\sum_{x\in\Lambda}d_{x}(\lambda)\theta_{x}$, where
$d_{x}(\lambda)=\textrm{dim}_{\mathbb{C}}(S_{\lambda})_{x}$. Thus we
have that
$\varphi(S_{\lambda})=\widetilde{\varphi}(S_{\lambda})=\sum_{x\in\Lambda}
d_{x}(\lambda)\widetilde{\varphi}(\theta_{x})=\\
\sum_{{p=0,1}\atop
x'\in\Lambda'}d_{x'}(\varphi(S_{\lambda}))\theta_{x'}'\epsilon^p$.
With this decomposition, we can easily get those $a_{k,p}$.

Now we have our main result in this section.

\noindent {\bf Theorem 3.6} For any $\lambda\in \Lambda^+$, we have
that $\varphi (S_{\lambda})=\sum_{{p=0,1}\atop k\in
\mathbb{N}}a_{k,p}V(k)\epsilon^p$, for some $a_{k,p}\in\mathcal
{A}$. Then we get
$$C_{rt}S_{\lambda}=\sum_{{p=0,1}\atop k\in \mathbb{N} }a_{k,p}C_{rt(srt)^k\omega^p}
\ \ \textrm{mod}\ \mathcal {H}_{c_0}.$$

\noindent {\bf Proof.} Since $\varphi (S_{\lambda})\in
R_{F_{c_2}}\otimes\mathcal {A}$, thus we can assume that $\varphi
(S_{\lambda})=\\ \sum_{{p=0,1}\atop k\in
\mathbb{N}}a_{k,p}V(k)\epsilon^p$, for some $a_{k,p}\in\mathcal
{A}$. By the definition of $\varphi_1,$ we get
$$\varphi_1^{-1}(\sum_{{p=0,1}\atop k\in
\mathbb{N}}a_{k,p}V(k)\epsilon^p)=\frac{1}{[2]^2}\sum_{{p=0,1}\atop
k\in \mathbb{N}}a_{k,p}\widehat{C}_{rt(srt)^{m}\omega^p}.$$ On the
other hand, we have that $$\varphi_1^{-1}(\sum_{{p=0,1}\atop k\in
\mathbb{N}}a_{k,p}V(k)\epsilon^p)=\varphi_1^{-1}\varphi
(S_{\lambda})=\varphi_2^{-1}(S_{\lambda})=\frac{1}{[2]^2}\widehat{C_{rt}S_{\lambda}}.$$
Thus we get that
$$C_{rt}S_{\lambda}=\sum_{{p=0,1}\atop k\in \mathbb{N} }a_{k,p}C_{rt(srt)^k\omega^p}
\ \ \textrm{mod}\ \mathcal {H}_{c_0}.$$ \hfill$\Box$\

\noindent {\bf Remark 3.7} \\
(1) In fact, the ideas in this section have an interpretation in
terms of based ring $\mathcal {J}$ introduced by Lusztig in [L2]. We
also need a conjecture on based ring $\mathcal {J}_{c}$ associated
with a two-sided cell $c$ of $W$ given by Lusztig. In case $W$ is
the affine Weyl group of type $\widetilde{B}_2$, this conjecture has
been proved by Xi. The details can be found in [X1, X3] or Chapter
11 in [X2].

\noindent (2) We also have some questions about the leading
coefficients of the Kazhdan-Lusztig polynomials for an extended
affine Weyl group of type $\widetilde{B_2}$. For those $y<w$
satisfying $a(y)=4$ and $a(w)=1$ or 2, we can only get part of the
leading coefficients. For example, we have computed that
$\mu(y,w)=0$ when $y$ and $w$ are both of the minimal length in
their double cosets $W_0yW_0$ and $W_0wW_0$, respectively. Then we
can know all $\mu(y,w)$ for those $y<w$ satisfying $l(w)-l(y)\leq3$.
These results will appear somewhere.

\noindent (3) The part contained in $\mathcal {H}_{c_0}$ of
$C_{rt}S_{\lambda}$ mainly depends on the solution to problem (2).

\subsection*{Acknowledgment} I must thank Professor N.Xi a lot for
his useful conversation and suggestions.

\providecommand{\bysame}{\leavevmode\hbox
to3em{\hrulefill}\thinspace}
\providecommand{\MR}{\relax\ifhmode\unskip\space\fi MR }
\providecommand{\MRhref}[2]{%
  \href{http://www.ams.org/mathscinet-getitem?mr=#1}{#2}
} \providecommand{\href}[2]{#2}

\end{document}